\documentclass[runningheads]{llncs}
\usepackage{amssymb}
\usepackage{amsmath}
\usepackage{bbold} 
\usepackage{color}
\usepackage{xcolor}
\usepackage{graphicx}
\usepackage{hyperref}
\def\alert{\mathbf}
\def\idmat{\bbbone}
\def\cmath{\color{black}}
\def\ctxt{\color{black}}
\newcommand{\AltG}[1]{\mathsf{A}_{#1}} 
\def\bmat{\begin{pmatrix}}
\def\emat{\end{pmatrix}}
\newcommand{\barket}[1]{\left|#1\right\rangle} 
\newcommand{\brabar}[1]{\left\langle#1\right|} 
\newcommand{\braket}[1]{\left\langle#1\right\rangle} 
\newcommand{\cabs}[1]{\left|#1\right|} 
\newcommand{\CyclG}[1]{\mathsf{C}_{#1}} 
\DeclareMathOperator*{\DirectSum}{\oplus}
\newcommand{\farg}[1]{\!\left(#1\right)} 
\def\Hspace{\mathcal{H}} 
\def\id{\mathbf{1}} 
\newcommand{\innerform}[3]{\braket{#1\cabs{#2}#3}} 
\newcommand{\IrrRep}[1]{\mathbf{#1}} 
\newcommand{\Math}[1]{$\cmath{}#1$} 
\newcommand{\MathEq}[1]{\begin{equation*}\cmath{#1}\end{equation*}}
\newcommand{\MathEqLab}[2]{\begin{equation}\cmath{#1}\label{#2}\end{equation}}
\newcommand{\Mthree}[9]{\bmat#1&#2&#3\\#4&#5&#6\\#7&#8&#9\emat} 
\newcommand{\Mtwo}[4]{\bmat#1&#2\\#3&#4\emat} 
\def\NF{\mathcal{F}} 
\newcommand{\ordset}[1]{\left[#1\right]} 
\newcommand{\PermRep}[1]{\mathbf{\underline{#1}}} 
\def\Q{\mathbb{Q}} 
\newcommand{\set}[1]{\left\{#1\right\}} 
\newcommand{\SymGN}[1]{\mathsf{S}_{#1}} 
\newcommand{\vect}[1]{\left(#1\right)} 
\DeclareMathOperator{\zeromat}{\mathbb{0}}
\begin{document}
\title{An Algorithm for Computing Invariant Projectors in Representations of Wreath Products}
\titlerunning{Invariant Projectors in Representations of Wreath Products}
\author{Vladimir V. Kornyak 
}
\institute{Laboratory of Information Technologies \\
           Joint Institute for Nuclear Research \\
           141980 Dubna, Russia \\
           \email{vkornyak@gmail.com}}
\authorrunning{Vladimir V. Kornyak}
\maketitle
\begin{abstract}
We describe an algorithm for computing the complete set of primitive orthogonal idempotents in the centralizer ring of the permutation representation of a wreath product.
This set of idempotents determines the decomposition of the representation into irreducible components.
In the formalism of quantum mechanics, these idempotents are projection operators into irreducible invariant subspaces of the Hilbert space of a multipartite quantum system.
The C implementation of the algorithm constructs irreducible decompositions of high-dimensional representations of wreath products.
Examples of computations are given.

\end{abstract}
\section{Introduction}\label{intro}
A typical description of a physical system usually includes a \emph{space} \Math{X}, on which a group of \emph{spatial symmetries} \Math{G=G\farg{X}} acts,
and a set of \emph{local states} \Math{V} with a group of \emph{local} (or \emph{internal}) \emph{symmetries} \Math{F=F\farg{V}}.
The sets \Math{X} and \Math{V} and the group \Math{F} can be thought of, respectively, as the \emph{base}, the \emph{typical fiber} and the \emph{structure group} of a \emph{fiber bundle}.
The \emph{sections} of the bundle are the set of functions from \Math{X} to \Math{V}, denoted by \Math{V^X}.
The set \Math{V^X} describes the whole states of the physical system.
A natural symmetry group that acts on \Math{V^X} and preserves the structure of the fiber bundle is the \emph{wreath product} \cite{Meldrum,JamesKerber} of \Math{F} and \Math{G}
\MathEq{\widetilde{W}=F\wr{}G\cong{}F^X\rtimes{}G\,.} 
The \emph{primitive}%
\footnote{Another canonical action of the wreath product, the \emph{imprimitive} action, acts on the fibers, i.e., on the set \Math{V\times{}X}.
We will not consider this action here.}
action of \Math{\widetilde{W}} on \Math{V^X} is described by the rule
\MathEqLab{v\farg{x}\vect{f\farg{x}, g}=v\farg{x{g^{-1}}}f\farg{x{g^{-1}}},}{wreathaction}
where \Math{v\in{}V^X}, \Math{f\in{}F^X}, \Math{g\in{}G}; 
the \emph{right-action} convention is used for all group actions.
\par
The wreath product is an important mathematical construct.
In particular, the \emph{universal embedding theorem} (also known as the \emph{Kaloujnine-Krasner theorem}) states that any extension of group \Math{A} by group \Math{B} is isomorphic to a subgroup of \Math{A\wr{}B}, i.e., the wreath product is a universal object containing all extensions.
Wreath products play an important role \cite{Cameron} in the influential \emph{O'Nan-Scott theorem}, which classifies \emph{maximal subgroups} of the \emph{symmetric} group, yet another universal object.
The wreath product \Math{\SymGN{m}\wr\SymGN{n}} is the automorphism group of the \emph{hypercubic graph} or \emph{Hamming scheme} \Math{H\farg{n, m}} in coding theory \cite{BannaiIto}.
\par
A quantum description of a physical system is obtained by introducing a Hilbert space \Math{\Hspace} spanned by the ``classical'' states of the system.
Accordingly, the action of the symmetry group on classical states goes over to the unitary representation of the group in the Hilbert space \Math{\Hspace}.
The next important and natural step in the study of the quantum system is the decomposition of the representation into irreducible components.
\par
Among the problems of quantum mechanics, the study of multipartite quantum systems is of particular interest because they manifest such phenomena as entanglement and non-local correlations.
In particular, the very possibility of quantum computing is based on these phenomena.
The Hilbert space of a multipartite system consisting of \Math{N} identical constituents has the form
\MathEqLab{\widetilde{\Hspace}=\Hspace^{\otimes{N}},}{expH}
where \Math{\Hspace} is the Hilbert space of a single constituent.
Assuming that the local group \Math{F} acts through a representation in \Math{\Hspace}, and the group \Math{G\leq\SymGN{N}} permutes the constituents, we come to the representation of the wreath product \Math{F\wr{}G} in the space \Math{\widetilde{\Hspace}}.
\par
In \cite{KornyakCASC18}, we proposed an algorithm for decomposing representations of finite groups into irreducible subrepresentations.
The algorithm computes a complete set of mutually orthogonal irreducible invariant projectors.
In fact, this is a special case of a general construction --- a \emph{Peirce decomposition} of a ring with respect to an orthogonal system of idempotents (see \cite{Jacobson,Rowen} for more details).
Invariant projection operators are important in problems of quantum physics, since they define invariant inner products in invariant subspaces of a Hilbert space.
Computer implementation of the approach has proved to be very effective in many cases.
For example, the program coped with many high dimensional representations of simple groups and their ``small'' extensions, presented in the \textsc{Atlas} \cite{atlas}, in the computationally difficult case of characteristic zero.
The algorithm in \cite{KornyakCASC18} uses polynomial algebra methods, which by their nature are algorithmically hard.
The number of polynomial variables is equal to the rank -- defined as the dimension of the centralizer ring -- of the representation to be split.
Representations of simple (or ``close'' to simple) groups usually have low ranks, and in such cases the algorithm works well.
However, wreath products, which contain all possible extensions, are far from simple groups and have high ranks.
Therefore, the approach from \cite{KornyakCASC18} is not applicable in the case of wreath products.
In this paper, we propose an algorithm for calculating irreducible invariant projectors in the representation of a wreath product in the form of tensor product polynomials with the projectors of a local group representation as variables.
We will consider here permutation representations -- the most fundamental type of representations.

\section{
Irreducible Invariant Projectors of Wreath Product}\label{algo}
Let \Math{X\cong\set{1,\ldots,N}} and  \Math{V\cong\set{1,\ldots,M}}.
This implies that \Math{G\farg{X}\leq{}\SymGN{N}} and \Math{F\farg{V}\leq{}\SymGN{M}}.
The functions \Math{v\in{}V^X} and \Math{f\in{}F^X} can be thought as arrays (ordered lists) \Math{\ordset{v_{1},\ldots,v_{N}}} and \Math{\ordset{f_{1},\ldots,f_{N}}}, respectively.
Accordingly, the wreath product element \Math{\widetilde{w}\in\widetilde{W}} can be written as the pair \Math{\bigl(\ordset{f_{1},\ldots,f_{N}}; g\bigr)}, where \Math{g\in{}G}.
The action \eqref{wreathaction} takes the form
\MathEq{\ordset{v_{1},\ldots,v_{N}}\xrightarrow{\bigl(\ordset{f_{1},\ldots,f_{N}};~ g\bigr)}\ordset{v_{1g^{-1}}f_{1g^{-1}},\ldots,v_{Ng^{-1}}f_{Ng^{-1}}}.}
The permutation representation \Math{\widetilde{P}} of the wreath product  is a representation of \Math{\widetilde{W}} by \Math{(0,1)}-matrices of the size  \Math{M^N\times{}M^N}  that have the form
\MathEqLab{\widetilde{P}\farg{\widetilde{w}}_{u,v}=\delta_{u\widetilde{w},v}, ~\text{\ctxt{}where}~ \widetilde{w}\in\widetilde{W};~u,v\in{}V^X;~ \delta~\text{\ctxt{}is the \emph{Kronecker delta}.}}{permrep}
We assume that the representation space is an \Math{M^N}-dimensional Hilbert space \Math{\widetilde{\Hspace}} over some abelian extension of the field \Math{\Q}.
This extension \Math{\NF} is a splitting field of the local group \Math{F}.
\subsection{Centralizer Ring of Wreath Product}
Let \MathEqLab{\widetilde{A}_1,\ldots,\widetilde{A}_{\widetilde{R}}}{basis}
be the basis elements of the centralizer ring of the representation \eqref{permrep}, \Math{\widetilde{R}} denotes the rank of the representation.
The basis elements \eqref{basis} are solutions of the system of equations (invariance condition)
\MathEqLab{\widetilde{P}\farg{\widetilde{w}^{-1}}\widetilde{A}\widetilde{P}\farg{\widetilde{w}}=\widetilde{A},~\widetilde{w}\in\widetilde{W}\,.}{invarcond}
A more detailed analysis of \eqref{invarcond} shows that the elements \eqref{basis} are in one-to-one correspondence with the orbits of \Math{\widetilde{W}} on the Cartesian square \Math{V^X\times{}V^X}. Such orbits are called \emph{orbitals}.
\par
Let us present the Cartesian square in the form \Math{\vect{V\times{}V}^X}, i.e.,  as the array \MathEqLab{\ordset{\vect{V\times{}V}_1,\ldots,\vect{V\times{}V}_N}.}{array}
To construct orbitals, consider the structure of the group \Math{\widetilde{W}=F\wr{}G} in more detail.
Its subgroup
\MathEqLab{\widetilde{F^X}=\vect{F^X;\,\id_G}\cong{}F^X}{basegroup}
is called the \emph{base group} of the wreath product.
The group \Math{\widetilde{F^X}} is a \emph{normal} (or \emph{invariant}) subgroup of \Math{\widetilde{W}}. This means that
\Math{\widetilde{w}^{-1}\widetilde{F^X}\widetilde{w}=\widetilde{F^X}} for any \Math{\widetilde{w}\in\widetilde{W}}, and this is denoted by \Math{\widetilde{F^X}\triangleleft\widetilde{W}}.
The subgroup \MathEqLab{\widetilde{G}=\vect{\id_{\!F}{\!\!}^X;\,G}\cong{}G}{G}
is a complement of \Math{\widetilde{F^X}} in \Math{\widetilde{W}}, i.e.,
\MathEq{\widetilde{W}=\widetilde{F^X}\cdot\widetilde{G} \text{\ctxt~~and~~} \widetilde{F^X}\cap\widetilde{G}=\id_{\widetilde{W}}\equiv\vect{\id_{\!F}{\!\!}^X;\,\id_G}.}
Thus, we can construct the orbits on the set \Math{\vect{V\times{}V}^X} acting first by the elements of the base group \eqref{basegroup}, and then by the elements of the complement \eqref{G}.
Further, we note that \Math{F^X} being the direct product of \Math{N} copies of \Math{F}, \Math{F^X=F_1\times\cdots\times{}F_N}, can be applied to the array \eqref{array} component wise independently. The action of the local group \Math{F} splits the set \Math{V\times{}V} into \Math{R} disjoint subsets
\MathEq{\Delta_1,\ldots,\Delta_R,}
which we will call \emph{local orbitals}.
Calculating local orbitals is a simple task, since the local group is exponentially smaller than the wreath product.
Let \Math{\overline{R}=\set{1,\ldots,R}} and \Math{\overline{R}^X} be the set of all mappings from \Math{X} into \Math{\overline{R}}.
We define the action of \Math{g\in{}G} on the mapping \Math{r\in\overline{R}^X} by \Math{rg=\ordset{r_{1g},\ldots,r_{Ng}}}.
This action decomposes the set \Math{\overline{R}^X} into orbits, and we can write the orbital of the wreath product as
\MathEq{\widetilde{\Delta}_r=\bigsqcup_{q\in{}rG}\Delta_{q_1}\times\cdots\times\Delta_{q_N}\,,} 
where \Math{rG} denotes the orbit of the mapping \Math{r}.
To translate from the language of sets to the language of matrices, we must replace the orbitals with the basis elements of the local centralizer ring, union by summation, and Cartesian products by the tensor products.
Thus, we obtain the expression for the basis element of the centralizer ring of the wreath product
\MathEqLab{\widetilde{A}_r=\sum_{q\in{}rG}A_{q_1}\otimes\cdots\otimes{}A_{q_N}\,,}{wreathA}
where \Math{A_1,\ldots,A_R} are basis elements of the local centralizer ring.
It is easy to show that the basis elements \eqref{wreathA} form a complete system, i.e.,
\MathEq{\sum_{i=1}^{\widetilde{R}}\widetilde{A}_{r^{(i)}}=\mathbb{J}_{M^N},}
where \Math{\mathbb{J}_{M^N}} is  the \Math{M^N\times{}M^N} \emph{all-ones matrix}, and \Math{r^{(i)}} denotes some numbering of the orbits of \Math{G} on \Math{\overline{R}^X}.
\subsection{Complete Set of Irreducible Invariant Projectors}
The complete set of irreducible invariant projectors is a subset of the centralizer ring, specified by the conditions of idempotency and mutual orthogonality.
A similar construction in ring theory is called a \emph{complete set of primitive orthogonal idempotents}.
An arbitrary ring with a complete set of orthogonal idempotents can be represented as a direct sum of indecomposable rings.
This is called a \emph{Peirce decomposition} \cite{Jacobson,Rowen}.
\par
Before constructing the complete set of primitive orthogonal idempotents for the centralizer ring of the permutation representation of the wreath product, we recall some properties of the tensor (Kronecker) product \cite{Steeb}:
\begin{enumerate}
	\item
\Math{\vect{A\otimes{}B}\otimes{}C=A\otimes\vect{B\otimes{}C}},
	\item
\Math{\vect{A+B}\otimes{}\vect{C+D}=A\otimes{}C+A\otimes{}D+B\otimes{}C+B\otimes{}D},
	\item
\Math{\vect{\alpha{}A}\otimes{}B=A\otimes\vect{\alpha{}B}=\alpha\vect{A\otimes{}B}},
	\item
\Math{\vect{A\otimes{}B}\vect{C\otimes{}D}=\vect{AC}\otimes\vect{BD},}
	\item
\Math{\vect{S\otimes{}T}^{-1}=S^{-1}\otimes{}T^{-1},}	
\end{enumerate}
where \Math{A}, \Math{B}, \Math{C} and \Math{D} are matrices, \Math{S} and \Math{T} are invertible matrices and \Math{\alpha} is a scalar.
It follows immediately from these properties that
\begin{enumerate}
	\item
if \Math{A}	and \Math{B} are both invariant, then \Math{A\otimes{}B} is invariant;
	\item
if \Math{A}	and \Math{B} are both idempotents, then \Math{A\otimes{}B} is idempotent;
	\item
if \Math{A'=S^{-1}AS}	and \Math{B'=T^{-1}BT}, then
\MathEq{A'\otimes{}B'=\vect{S\otimes{}T}^{-1}\vect{A\otimes{}B}\vect{S\otimes{}T}\equiv\vect{A\otimes{}B}'.}
This relation means that we can freely change the bases in the factors of the tensor product to more convenient ones.
\end{enumerate}
Using the above relations, their consequences and some additional technical considerations allows us to come to the main result of this section.
\par
Let \Math{B_1,\ldots,B_K} be the complete set of irreducible invariant projectors of the permutation representation of the local group \Math{F}.
Let \Math{\overline{K}=\set{1,\ldots,K}} and \Math{\overline{K}^X} be the set of all mappings from \Math{X} into \Math{\overline{K}}.
The action of \Math{g\in{}G} on the mapping \Math{k\in\overline{K}^X} is defined as \Math{kg=\ordset{k_{1g},\ldots,k_{Ng}}}.
Then we have
\begin{proposition}
An irreducible invariant projector in the permutation representation of the wreath product takes the form
\MathEqLab{\widetilde{B}_k=\sum_{\ell\in{}kG}B_{\ell_1}\otimes\cdots\otimes{}B_{\ell_N}\,,}{wreathB}
where \Math{kG} denotes the \Math{G}-orbit of the mapping \Math{k} on the set \Math{\overline{K}^X}.
\end{proposition}
The easily verifiable completeness condition
\Math{\displaystyle\sum_{i=1}^{\widetilde{K}}\widetilde{B}_{k^{(i)}}=\idmat_{M^N}}
holds.
Here \Math{\widetilde{K}} is the number of irreducible components of the wreath product representation,
\Math{\idmat_{M^N}} is the identity matrix in the representation space, \Math{k^{(i)}} denotes some numbering of the orbits of \Math{G} on \Math{\overline{K}^X}.
\par
To calculate the basis elements \eqref{wreathA} of the centralizer ring and irreducible invariant projectors \eqref{wreathB}, we wrote a program in C.
The input data for the program are the generators of the spatial and local groups, and the complete set of irreducible invariant projectors of the local group (obtained, for example, using the program described in \cite{KornyakCASC18}).
\section{Calculation Example: \Math{\SymGN{4}\farg{\text{\textit{octahedron}}}\wr\AltG{5}\farg{\text{\textit{icosahedron}}}}}\label{imple}
We give here the calculation of the centralizer ring and invariant projectors for the permutation representation of the wreath product of the rotational symmetry groups of the octahedron and icosahedron. This representation has dimension \Math{M^N={2\,176\,782\,336}} and rank \Math{\widetilde{R}={9099}}\,.
\subsection{Space Group}
\begin{figure}
\centering
\includegraphics[height=0.4\textwidth]{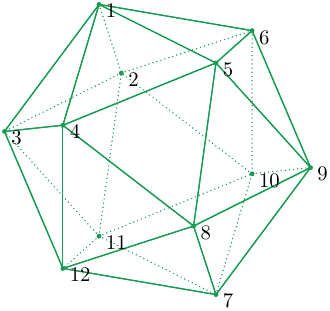}
\caption{Icosahedron.} \label{Icosahedron}
\end{figure}

In our example, the space \Math{X} is represented by the icosahedron, see Fig. \ref{Icosahedron}.
The full symmetry group of the icosahedron is the product \Math{\AltG{5}\times\CyclG{2}}.
As a group of spatial symmetries we take the orientation-preserving factor \Math{\AltG{5}},
which describes the \emph{rotational} (or \emph{chiral}) \emph{symmetries} of the icosahedron.
The order of \Math{\AltG{5}} is equal to \Math{60}.
The points of the space \Math{X} are the vertices of the icosahedron: \Math{X\cong\set{1,\ldots,12}}.
The space symmetry group \Math{G\farg{X}} can be generated by two permutations.
For example, for the vertex numbering as in Fig. \ref{Icosahedron}, we may use the following presentation
\MathEq{G\farg{X}=\bigl\langle(1,7)(2,8)(3,12)(4,11)(5,10)(6,9),\,(1,2,11,12,4)(5,6,10,7,8)\bigr\rangle\,\cong\,\AltG{5}.}
\subsection{Local Group}
\begin{figure}
\centering
\includegraphics[height=0.4\textwidth]{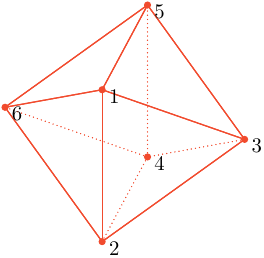}
\caption{Octahedron.} \label{Octahedron}
\end{figure}

The  set of local states is represented by the vertices of the octahedron (Fig. \ref{Octahedron}): \Math{V\cong\set{1,\ldots,6}}.
The full octahedral symmetry group is isomorphic to \Math{\SymGN{4}\times\CyclG{2}}.
Note that the full octahedral group itself is a wreath product, namely, \Math{\SymGN{2}\wr\SymGN{3}}.
The group of rotational symmetries of the octahedron is \Math{\SymGN{4}}.
The order of \Math{\SymGN{4}} is equal to \Math{24}.
For the vertex numbering of Fig. \ref{Octahedron}, the local symmetry group \Math{F\farg{V}} has the following presentation by two generators
\MathEq{F\farg{V}=\bigl\langle\vect{1,3,5}\vect{2,4,6},\,\vect{1,2,4,5}\bigr\rangle\,\cong\,\SymGN{4}.}
The six-dimensional permutation representation \Math{\PermRep{6}} of \Math{F\farg{V}} has rank \Math{3},
and the basis of the centralizer ring is
\MathEqLab{A_1=\idmat_6,~A_2=\Mtwo{\zeromat_3}{\idmat_3}{\idmat_3}{\zeromat_3},~A_3=\Mtwo{Y}{Y}{Y}{Y}, \text{\ctxt~where~} Y=\Mthree{0}{1}{1}{1}{0}{1}{1}{1}{0}.}{localA}
The irreducible decomposition of this representation has the form \Math{\PermRep{6}=\IrrRep{1}\DirectSum\IrrRep{2}\DirectSum\IrrRep{3}}.
The complete set of irreducible invariant projectors for representation \Math{\PermRep{6}}, expressed in the basis \eqref{localA} of the centralizer ring, is as follows
\MathEqLab{B_{\IrrRep{1}}=\frac{1}{6}\vect{A_1+A_2+A_3}, B_{\IrrRep{2}}=\frac{1}{3}\vect{A_1+A_2-\frac{1}{2}A_3}, B_{\IrrRep{3}}=\frac{1}{2}\vect{A_1-A_2}.}{localB}

\subsection{Applying C Program}
Below we present the output of the program, accompanying parts of the output with comments.
An example with a larger dimension of representation is given in Appendix \ref{appendix}.
Both calculations were performed on PC with 3.30GHz CPU and 16GB RAM.
\subsubsection{Part 1.}
The program reads the generators of the spatial and local groups and prints related information:
`\texttt{Name}' is the name of the file containing the generators,
`\texttt{Number of points}' is the dimension of the permutation representation, etc.
~\\[-20pt]
\begin{verbatim}
Space G(X) group:
  Name = "A5_on_icosahedron"
  Number of points = 12
  Comment = "Action of A_5 on 12 vertices of icosahedron"
  Size = "60"
  Number of generators = 2
Local F(V) group:
  Name = "S4_on_octahedron"
  Number of points = 6
  Comment = "Action of S_4 on 6 vertices of octahedron"
  Size = "24"
  Number of generators = 2
\end{verbatim}
\subsubsection{Part 2.}
The program constructs the wreath product generators from the generators of the constituent groups and creates a file for the constructed generators.
This part of the calculation is optional and can be disabled.
\begin{verbatim}
Whole F(V).wr.G(X) group:
  Name = "S4_on_octahedron_wr_A5_on_icosahedron"
  Number of points V^X = 2176782336
  Size = "24^12*60"
  Number of generators = 4
\end{verbatim}
\subsubsection{Part 3.}
The program prints a common header (in \LaTeX):\\
\textbf{Wreath product} \Math{\SymGN{4}\farg{\text{\textit{octahedron}}}\wr\AltG{5}\farg{\text{\textit{icosahedron}}}}\\
Representation dimension: \Math{\alert{2\,176\,782\,336}}
\subsubsection{Part 4.} The program computes the basis elements of the centralizer ring of the wreath product.
This optional computation can be skipped.
If it is enabled, the following items are  calculated:
\begin{enumerate}
	\item
\emph{Rank} of the representation (dimension of the centralizer ring).
	\item
The set of \emph{suborbit lengths}.
For a permutation group, suborbits are the orbits of the stabilizer of a single point of the permutation domain \cite{Cameron}.
The sizes of orbitals are equal to the suborbit lengths multiplied by the representation dimension.
The sum of suborbit length is equal to the  representation dimension ---
we use this fact to verify calculation by displaying `\texttt{Checksum}'.
The multiplicities of the suborbit lengths are shown as superscripts.
	\item
The basis elements of the centralizer ring for the wreath product representation.
We display here only a few elements of the large array.
The expressions for the basis elements are tensor product polynomials in the matrices \Math{A_1, A_2} and \Math{A_3}, shown explicitly in \eqref{localA}.
\end{enumerate}
{\normalsize{}
\texttt{Rank:} \Math{\alert{9099}}\\
\texttt{Number of different suborbit lengths:} \Math{37}\\
\texttt{Wreath suborbit lengths:}\\[3pt]
{\Math{1^{2}, 6^{2}, 10^{2}, 12^{6}, 15^{2}, 20^{6}, 30^{24}, 48^{8}, 60^{52}, 96^{2}, 192, 240^{408}, 480^{94}, 960^{1079}, 1280^{16},}\\
{\Math{3840^{1876}, 7680^{114}, 12288^{8}, 15360^{\mathbf{\underline{2054}}}, 40960^{4}, 49152^{8}, 61440^{1688}, 81920^{10}, 122880^{76},}\\
\Math{196608^{2}, 245760^{942}, 983040^{424}, 1966080^{27}, 3932160^{118}, 5242880^{4}, 6291456^{2},}\\
\Math{12582912, 15728640^{28}, 16777216, 31457280^{4}, 50331648^{2}, \mathbf{\underline{62914560}}^{2}}\\
\texttt{Checksum =} \Math{2176782336}  \texttt{Maximum multiplicity =} \Math{2054}\\[2pt]
\texttt{Wreath invariant basis forms:}
\begin{align*}
\cmath\widetilde{A}_{1}=&\cmath{}A_{1}^{\otimes12}\\
\cmath\widetilde{A}_{2}=&\cmath{}A_{2}^{\otimes12}\\
\cmath\widetilde{A}_{3}=&\cmath{}A_{1}^{\otimes{}5}\otimes{}A_{2}\otimes{}A_{1}^{\otimes{}5}\otimes{}A_{2} + A_{1}^{\otimes{}4}\otimes{}A_{2}\otimes{}A_{1}^{\otimes{}5}\otimes{}A_{2}\otimes{}A_{1} + A_{1}^{\otimes{}3}\otimes{}A_{2}\otimes{}A_{1}^{\otimes{}5}\otimes{}A_{2}\otimes{}A_{1}^{\otimes{}2}\\
+&A_{1}^{\otimes{}2}\otimes{}A_{2}\otimes{}A_{1}^{\otimes{}5}\otimes{}A_{2}\otimes{}A_{1}^{\otimes{}3} + A_{1}\otimes{}A_{2}\otimes{}A_{1}^{\otimes{}5}\otimes{}A_{2}\otimes{}A_{1}^{\otimes{}4} + A_{2}\otimes{}A_{1}^{\otimes{}5}\otimes{}A_{2}\otimes{}A_{1}^{\otimes{}5}\\
\cmath\vdots~&
\end{align*}
\subsubsection{Part 5.}
This main part of the program computes the complete set of irreducible invariant projectors \eqref{wreathB}.
The output below contains
\begin{enumerate}
	\item
Information on whether the irreducible decomposition is multiplicity free or contains subrepresentations with nontrivial multiplicity.	
	\item
The total number of the irreducible invariant projectors.	
	\item
The number of different dimensions.
\item
Complete list of irreducible dimensions.
The superscripts represent the numbers of irreducible projectors of respective dimensions.
The total sum of dimensions must be equal to the dimension of the wreath product representation,
and this is verified by the direct calculation of `\texttt{Checksum}'.
\item
A few expressions for irreducible invariant projectors, which are tensor product polynomials in matrices of local projectors \eqref{localB}.
\end{enumerate}
\texttt{Wreath product decomposition is multiplicity free}\\
\texttt{Number of irreducible components:} \Math{9099}\\
\texttt{Number of different dimensions:} \Math{125}\\
\texttt{Irreducible dimensions:}\\
{\Math{1, 24^{2}, 36, 54, 72, 120^{2}, 160^{2}, 240, 270^{2}, 360^{2}, 384, 480^{9}, 540^{2}, 640^{2}, 720^{11}, 768^{2}, 960^{5},}\\
\Math{1080^{26}, 1152^{2}, 1215, 1280^{2}, 1440^{33}, 1536, 1620^{3}, 1920^{21}, 2160^{42},2304^{2}, 2430^{6}, 2880^{66},}\\
\Math{2916, 3240^{33}, 3456, 3840^{11}, 4096, 4320^{168}, 4860^{5}, 5760^{92}, 5832^{2}, 6144, 6480^{132}, 7290^{2},}\\
\Math{7680^{19}, 8640^{216}, 8748^{2}, 9720^{96}, 10240^{2}, 11520^{92}, 11664, 12960^{306}, 14580^{15}, 15360^{5}, }\\
\Math{17280^{307}, 17496^{2}, 19440^{218}, 21870^{8}, 23040^{66}, 24576, 25920^{462}, 26244, 29160^{92}, 30720^{5},}\\
\Math{34560^{222}, 36864, 38880^{504}, 43740^{10}, 46080^{33}, 51840^{462}, 55296, 58320^{277}, 69120^{147}, 73728,}\\
\Math{77760^{\alert{\underline{556}}}, 87480^{122}, 92160^{11}, 93312, 98415, 103680^{306}, 116640^{468}, 131220^{13}, 138240^{42},}\\
\Math{155520^{492}, 174960^{216}, 184320^{2}, 186624^{2}, 196830^{6}, 207360^{132}, 233280^{462}, 262440^{66}, 276480^{15},}\\
\Math{279936^{2}, 311040^{217}, 349920^{336}, 354294, 373248, 393660^{7}, 414720^{33}, 466560^{279}, 524880^{132},}\\
\Math{531441, 559872^{2}, 622080^{72}, 699840^{216}, 787320^{48}, 829440^{3}, 839808, 933120^{94}, 1049760^{132},}\\
\Math{1180980^{3}, 1244160^{5}, 1399680^{100}, 1417176^{2}, 1574640^{43}, 1771470^{2}, 1866240^{13}, 2099520^{66}, }\\
\Math{2125764,2361960^{11}, 2799360^{10}, 3149280^{41}, 4199040^{13}, 4251528, 4723920^{11}, 6298560^{5},}\\
\Math{7085880^{4}, \alert{\underline{9447840}}^{3}}\\[2pt]
\texttt{Checksum =} \Math{2176782336}  \texttt{Maximum number of equal dimensions =} \Math{556}\\[5pt]
\texttt{Wreath irreducible projectors:}
\begin{align*}
\cmath\widetilde{B}_{1}=&\cmath{}B_{1}^{\otimes12}\\
\cmath\widetilde{B}_{2}=&B_{1}^{\otimes{}11}\otimes{}B_{2} + B_{1}^{\otimes{}10}\otimes{}B_{2}\otimes{}B_{1} + B_{1}^{\otimes{}9}\otimes{}B_{2}\otimes{}B_{1}^{\otimes{}2} + B_{1}^{\otimes{}8}\otimes{}B_{2}\otimes{}B_{1}^{\otimes{}3}\\
+&  B_{1}^{\otimes{}7}\otimes{}B_{2}\otimes{}B_{1}^{\otimes{}4}+  B_{1}^{\otimes{}6}\otimes{}B_{2}\otimes{}B_{1}^{\otimes{}5} + B_{1}^{\otimes{}5}\otimes{}B_{2}\otimes{}B_{1}^{\otimes{}6} + B_{1}^{\otimes{}4}\otimes{}B_{2}\otimes{}B_{1}^{\otimes{}7}\\
+& B_{1}^{\otimes{}3}\otimes{}B_{2}\otimes{}B_{1}^{\otimes{}8} + B_{1}^{\otimes{}2}\otimes{}B_{2}\otimes{}B_{1}^{\otimes{}9}
 + B_{1}\otimes{}B_{2}\otimes{}B_{1}^{\otimes{}10} + B_{2}\otimes{}B_{1}^{\otimes{}11}\\
\widetilde{B}_{3}=&B_{1}^{\otimes{}5}\otimes{}B_{2}\otimes{}B_{1}^{\otimes{}5}\otimes{}B_{2} + B_{1}^{\otimes{}4}\otimes{}B_{2}\otimes{}B_{1}^{\otimes{}5}\otimes{}B_{2}\otimes{}B_{1} + B_{1}^{\otimes{}3}\otimes{}B_{2}\otimes{}B_{1}^{\otimes{}5}\otimes{}B_{2}\otimes{}B_{1}^{\otimes{}2}\\
+& B_{1}^{\otimes{}2}\otimes{}B_{2}\otimes{}B_{1}^{\otimes{}5}\otimes{}B_{2}\otimes{}B_{1}^{\otimes{}3} + B_{1}\otimes{}B_{2}\otimes{}B_{1}^{\otimes{}5}\otimes{}B_{2}\otimes{}B_{1}^{\otimes{}4} + B_{2}\otimes{}B_{1}^{\otimes{}5}\otimes{}B_{2}\otimes{}B_{1}^{\otimes{}5}\\
\cmath\vdots~&\\
\end{align*}
\texttt{Time: 0.89  sec}\\
\texttt{Maximum number of tensor monomials: 1062882}

\section{Application Remarks}
\label{concl}
One of the main goals of the work was to develop a tool for the study of models of multipartite quantum systems.
The projection operators obtained by the program are supposed to be used to calculate quantum correlations in such models.
These operators are matrices of huge dimension (for example, about ten trillion for projectors in Appendix A).
Obviously, the explicit calculation of such matrices is impossible.
However, the expression of projectors for wreath products in the form of tensor polynomials makes it possible to reduce the computation of quantum correlations to a sequence of computations with small matrices of local projectors.
To demonstrate this, recall that the computation of quantum correlations is ultimately reduced to the computation of scalar products.
Let \Math{\displaystyle\widetilde{\Phi}=\sum_{m\in{}V^{X}}\varphi_{m}\barket{e_{m_1}}\otimes\cdots\otimes\barket{e_{m_N}}\in\widetilde{\Hspace}}  and
\Math{\displaystyle\widetilde{\Psi}=\sum_{m\in{}V^{X}}\psi_{n}\barket{e_{n_1}}\otimes\cdots\otimes\barket{e_{n_N}}\in\widetilde{\Hspace}} be vectors
of the Hilbert space \eqref{expH}, where \Math{\barket{e_{1}},\ldots,\barket{e_{M}}} is a basis in the local Hilbert space \Math{\Hspace},
\Math{\varphi_{m}} and \Math{\psi_{n}} are arbitrary scalars from the base field \Math{\NF}.
Then we can calculate the scalar product of these vectors in the invariant subspace defined by the projector \eqref{wreathB} as follows
{\cmath\begin{multline*}
\innerform{\widetilde{\Phi}}{\widetilde{B}_k}{\widetilde{\Psi}}\\
=\sum_{\substack{m\in{}V^X\\
n\in{}V^X\\
\ell\in{}kG}
}\overline{\varphi}_m\psi_n\brabar{e_{m_1}}\otimes\cdots\otimes\brabar{e_{m_N}}
B_{\ell_1}\otimes\cdots\otimes{}B_{\ell_N}\barket{e_{n_1}}\otimes\cdots\otimes\barket{e_{n_N}}\\
=\sum_{\substack{m\in{}V^X\\
n\in{}V^X\\
\ell\in{}kG}
}\overline{\varphi}_m\psi_n\brabar{e_{m_1}}B_{\ell_1}\barket{e_{n_1}}\cdots\brabar{e_{m_N}}B_{\ell_N}\barket{e_{n_N}}\,.
\end{multline*}}

\subsection*{Acknowledgments}
I am grateful to Yu.A. Blinkov for help in preparing the article,  
V.P. Gerdt for fruitful discussions and valuable advice, 
and A. R\"omer for detecting a bug in the computer program..
\appendix
\section{Computing Invariant Projectors for \Math{\AltG{5}\farg{\text{\textit{icosahedron}}}\wr\AltG{5}\farg{\text{\textit{icosahedron}}}}}\label{appendix}

The representation \Math{\PermRep{12}} of \Math{\AltG{5}} on the  icosahedron  has rank \Math{4},
and the basis of the centralizer ring is
\MathEqLab{A_{1}=\idmat_{12},~ A_{2}=\Mtwo{\zeromat_{6}}{\idmat_{6}}{\idmat_{6}}{\zeromat_{6}},~A_{3}=\Mtwo{Y}{Z}{Z}{Y},~A_{4}=\Mtwo{Z}{Y}{Y}{Z},}{localA5}
where~ \Math{
Y=
\bmat
0&1&1&1&1&1\\
1&0&1&0&0&1\\
1&1&0&1&0&0\\
1&0&1&0&1&0\\
1&0&0&1&0&1\\
1&1&0&0&1&0\\
\emat,~
Z=
\bmat
0&0&0&0&0&0\\
0&0&0&1&1&0\\
0&0&0&0&1&1\\
0&1&0&0&0&1\\
0&1&1&0&0&0\\
0&0&1&1&0&0\\[2pt]
\emat.
}
The irreducible decomposition of this representation has the form \Math{\PermRep{12}=\IrrRep{1}\DirectSum\IrrRep{3}\DirectSum\IrrRep{3'}\DirectSum\IrrRep{5}}.
The complete set of irreducible invariant projectors for representation \Math{\PermRep{12}}, expressed in the basis \eqref{localA5} of the centralizer ring, is as follows
{\cmath	
\begin{align*}
	B_{\IrrRep{1}} =& \frac{1}{12}\vect{A_{1}+A_{2}+A_{3}+A_{4}},\\
	B_{\IrrRep{3}} =& \frac{1}{4}\vect{A_{1}-A_{2}+\frac{1}{\sqrt{5}}A_{3}-\frac{1}{\sqrt{5}}A_{4}},\\
	B_{\IrrRep{3'}} =& \frac{1}{4}\vect{A_{1}-A_{2}-\frac{1}{\sqrt{5}}A_{3}+\frac{1}{\sqrt{5}}A_{4}},\\
	B_{\IrrRep{5}} =& \frac{1}{12}\vect{5A_{1}+5A_{2}-A_{3}-A_{4}.}
\end{align*}
}

\begin{verbatim}
Local F(V) group:
  Name = "A5_on_icosahedron"
  Number of points = 12
  Comment = "Action of A_5 on 12 vertices of icosahedron"
  Number of generators = 2
\end{verbatim}
\begin{verbatim}
Space G(X) group:
  Name = "A5_on_icosahedron"
  Number of points = 12
  Comment = "Action of A_5 on 12 vertices of icosahedron"
  Number of generators = 2
\end{verbatim}
\begin{verbatim}
Whole group F(V).wr.G(X)
  Number of points V^X = 8916100448256
\end{verbatim}
\textbf{Wreath product} \Math{\AltG{5}\farg{\text{\textit{icosahedron}}}\wr\AltG{5}\farg{\text{\textit{icosahedron}}}}\\
\texttt{Representation dimension:} \Math{\alert{8\,916\,100\,448\,256}}\\
\texttt{Rank:} \Math{\alert{3\,875\,157}}\\
\texttt{Wreath product decomposition is multiplicity free}\\
\texttt{Number of irreducible components:} \Math{280\,832}\\
\texttt{Number of different dimensions:} \Math{145}\\
\texttt{Irreducible dimensions:}\\
\Math{1,36^{2},54^{2},60,108,150,180^{2},270^{4},540^{6},750^{2},900^{4},1215^{2},1620^{28},2430^{27},2500^{2},2700^{44},2916^{2},}\\
\Math{4500^{22},4860^{118},6750^{30},7290^{4},7500^{3},8100^{264},8748^{8},9375,13500^{183},14580^{436},18750^{6},21870^{76},}\\
\Math{22500^{66},24300^{1056},26244^{8},37500^{6},40500^{1056},43740^{950},60750^{120},67500^{540},72900^{2958},98415^{4},}\\
\Math{112500^{136},121500^{3636},131220^{1688},156250^{2},168750^{60},187500^{15},196830^{114},202500^{2460},218700^{5912},}\\
\Math{312500^{2},337500^{898},354294^{2},364500^{8870},393660^{2070},468750^{8},531441^{2},546750^{240},562500^{188},}\\
\Math{607500^{7392},656100^{8448},708588,759375^{4},937500^{11},1012500^{3696},1093500^{14664},1180980^{1872}\!,1518750^{174}\!,}\\
\Math{1687500^{1108},1771470^{94},1822500^{14808},1968300^{8448},2125764^{8},2812500^{184},3037500^{9152},3188646^{2},}\\
\Math{3280500^{16896},3542940^{1083},4218750^{60},4687500^{13},4920750^{240},5062500^{3696},5314410^{2},5467500^{\alert{\underline{19704}}},}\\
\Math{5859375,5904900^{5632},6377292^{6},7971615^{2},8437500^{906},8857350^{2},9112500^{14786},9841500^{12552},}\\
\Math{10628820^{422},11718750^{6},13668750^{240},14062500^{132},15187500^{7392},15943230^{24},16402500^{16896},}\\
\Math{17714700^{2253},23437500^{5},25312500^{2460},27337500^{14672},29524500^{5632},31886460^{52},37968750^{120},}\\
\Math{39062500^{2},42187500^{528},44286750^{94},45562500^{8874},49207500^{8464},53144100^{408},58593750,61509375^{4},}\\
\Math{70312500^{66},75937500^{3636},82012500^{8456},88573500^{1079},105468750^{30},113906250^{4},117187500^{3},}\\
\Math{123018750^{114},126562500^{1056},136687500^{5920},147622500^{1872},210937500^{183},227812500^{2968},244140625,}\\
\Math{246037500^{2054},292968750^{2},341718750^{76},351562500^{24},379687500^{1056},410062500^{1688},474609375^{2},}\\
\Math{527343750^{2},585937500,632812500^{264},683437500^{942},949218750^{27},1054687500^{49},1139062500^{422},}\\
\Math{1757812500^{6},1898437500^{118},2636718750^{4},3164062500^{28},\alert{\underline{5273437500}}^{2}}\\[2pt]
\texttt{Checksum =} \Math{8916100448256}  \texttt{Maximum number of equal dimensions =} \Math{19704}\\[5pt]
\texttt{Wreath irreducible projectors:}
{\cmath\begin{align*}
\widetilde{B}_{1}=&B_{\IrrRep{1}}^{\otimes12}\\
\widetilde{B}_{2}=&B_{\IrrRep{1}}^{\otimes11}\otimes{}B_{\IrrRep{3}}+B_{\IrrRep{1}}^{\otimes10}\otimes{}B_{\IrrRep{3}}\otimes{}B_{\IrrRep{1}}+B_{\IrrRep{1}}^{\otimes9}\otimes{}B_{\IrrRep{3}}\otimes{}B_{\IrrRep{1}}^{\otimes2}+B_{\IrrRep{1}}^{\otimes8}\otimes{}B_{\IrrRep{3}}\otimes{}B_{\IrrRep{1}}^{\otimes3}\\
+&B_{\IrrRep{1}}^{\otimes7}\otimes{}B_{\IrrRep{3}}\otimes{}B_{\IrrRep{1}}^{\otimes4}+B_{\IrrRep{1}}^{\otimes6}\otimes{}B_{\IrrRep{3}}\otimes{}B_{\IrrRep{1}}^{\otimes5}+B_{\IrrRep{1}}^{\otimes5}\otimes{}B_{\IrrRep{3}}\otimes{}B_{\IrrRep{1}}^{\otimes6}+B_{\IrrRep{1}}^{\otimes4}\otimes{}B_{\IrrRep{3}}\otimes{}B_{\IrrRep{1}}^{\otimes7}\\
+&B_{\IrrRep{1}}^{\otimes3}\otimes{}B_{\IrrRep{3}}\otimes{}B_{\IrrRep{1}}^{\otimes8}+B_{\IrrRep{1}}^{\otimes2}\otimes{}B_{\IrrRep{3}}\otimes{}B_{\IrrRep{1}}^{\otimes9}+B_{\IrrRep{1}}\otimes{}B_{\IrrRep{3}}\otimes{}B_{\IrrRep{1}}^{\otimes10}+B_{\IrrRep{3}}\otimes{}B_{\IrrRep{1}}^{\otimes11}\\
\widetilde{B}_{3}=&B_{\IrrRep{1}}^{\otimes11}\otimes{}B_{\IrrRep{3'}}+B_{\IrrRep{1}}^{\otimes10}\otimes{}B_{\IrrRep{3'}}\otimes{}B_{\IrrRep{1}}+B_{\IrrRep{1}}^{\otimes9}\otimes{}B_{\IrrRep{3'}}\otimes{}B_{\IrrRep{1}}^{\otimes2}+B_{\IrrRep{1}}^{\otimes8}\otimes{}B_{\IrrRep{3'}}\otimes{}B_{\IrrRep{1}}^{\otimes3}\\
+&B_{\IrrRep{1}}^{\otimes7}\otimes{}B_{\IrrRep{3'}}\otimes{}B_{\IrrRep{1}}^{\otimes4}+B_{\IrrRep{1}}^{\otimes6}\otimes{}B_{\IrrRep{3'}}\otimes{}B_{\IrrRep{1}}^{\otimes5}+B_{\IrrRep{1}}^{\otimes5}\otimes{}B_{\IrrRep{3'}}\otimes{}B_{\IrrRep{1}}^{\otimes6}+B_{\IrrRep{1}}^{\otimes4}\otimes{}B_{\IrrRep{3'}}\otimes{}B_{\IrrRep{1}}^{\otimes7}\\
+&B_{\IrrRep{1}}^{\otimes3}\otimes{}B_{\IrrRep{3'}}\otimes{}B_{\IrrRep{1}}^{\otimes8}+B_{\IrrRep{1}}^{\otimes2}\otimes{}B_{\IrrRep{3'}}\otimes{}B_{\IrrRep{1}}^{\otimes9}+B_{\IrrRep{1}}\otimes{}B_{\IrrRep{3'}}\otimes{}B_{\IrrRep{1}}^{\otimes10}+B_{\IrrRep{3'}}\otimes{}B_{\IrrRep{1}}^{\otimes11}\\
\vdots~&
\end{align*}}
\texttt{Time: 7.79 sec}\\
\texttt{Maximum number of tensor monomials: 16777216}

\end{document}